\documentclass[12pt]{article}
\usepackage{hyperref}
\usepackage{cmap}                        
\usepackage[cp1251]{inputenc}            
\usepackage[english]{babel}
\usepackage[left=1in,right=1in,top=1in,bottom=1.5in]{geometry} 
\usepackage{amssymb,amsmath, amsthm, amscd,ifthen}
\usepackage{graphicx}
\usepackage[svgnames]{xcolor}
\usepackage{tikz}
\usepackage{bm}
\usetikzlibrary{patterns,decorations.text,positioning,arrows,shapes,decorations.markings}
\tikzstyle{vertex}=[circle,draw=black,fill=black,inner sep=0,minimum size=0.2cm,text=white,font=\footnotesize]

\date{}

\title{Tilings of the plane with unit area triangles of bounded diameter}

\author{Andrey Kupavskii\thanks{University of Birmingham and MIPT, Moscow.  E-mail: {\tt kupavskii@ya.ru}.} \and J\'anos Pach\thanks{EPFL, Lausanne and R\'enyi Institute, Budapest. Supported by Swiss National Science Foundation Grants 200020-162884 and 200021-165977. E-mail: {\tt pach@cims.nyu.edu}.} \and G\'abor Tardos\thanks{R\'enyi Institute, Budapest. Supported by the Cryptography ``Lend\"ulet'' project of the Hungarian Academy
of Sciences and by the National Research, Development and Innovation Office, NKFIH, projects K-116769 and SNN-117879.}}

\date{}

\begin{document}
\maketitle
\centerline{\emph{\small Dedicated to the 70th birthdays of Ted Bisztriczky, G\'abor Fejes T\'oth, and Endre Makai}}
\begin{abstract}\noindent
There exist tilings of the plane with pairwise noncongruent triangles of equal area and bounded perimeter. Analogously, there exist tilings with triangles of equal perimeter, the areas of which are bounded from below by a positive constant. This solves a problem of Nandakumar.
\end{abstract}

\section{Introduction}
Nandakumar's questions on tilings~\cite{Na06, Na14} inspired a lot of research in geometry and topology. In particular, he and Ramana Rao~\cite{NaR12} conjectured that for every natural number $n$, any plane convex body can be partitioned into $n$ convex pieces of equal area and perimeter. After some preliminary results~\cite{BBS10} indicating that the problem is closely related to questions in equivariant topology, the conjecture was settled in the affirmative by Karasev, Hubard, and Aronov \cite{KaHA14} and independently by Blagojevi\'c and Ziegler~\cite{BlZ14}, in the special case where $n$ is a prime power. The problem remained open, e.g., for $n=6$. See also~\cite{So17, Zi15}.
\smallskip

Nandakumar also asked whether it possible to tile the plane with pairwise noncongruent triangles of the same area and perimeter. The answer turned out to be negative~\cite{KPT17}. In fact, we proved the stronger statement that every tiling of the plane with unit perimeter triangles such that the area of each triangle is at least $\varepsilon>0$, contains two triangles that share a side. If the areas of these two triangles are the same, they must be congruent.
\smallskip

This result left open the possibility that there exist tilings with pairwise noncongruent unit perimeter triangles whose areas are bounded from below by a constant $\varepsilon>0$. The ``dual'' question was also asked by Nandakumar: do there exist tilings of the plane with pairwise noncongruent unit area triangles with bounded perimeter? The aim of the present note is to show that the answer to both of these questions is positive.
\medskip

\noindent{\bf Theorem 1.} {\em There exist tilings of the plane with pairwise noncongruent triangles of unit area, the perimeters of which are bounded from above by a constant.}
\medskip

\noindent{\bf Theorem 2.} {\em There exist tilings of the plane with pairwise noncongruent triangles of unit perimeter, the areas of which are bounded from above by a constant.}

\medskip

Constructions very similar to ours have been found independently by Frettl\"oh~\cite{Fr16}. However, his proof has gaps. In particular, he argued that with one degree of freedom (that is, having uncountably many possible choices), one can avoid a countable set of ``coincidences''. While such a statement is true in many settings, it certainly does not hold in general. As our current proof shows, making such an argument work requires nontrivial effort.
\smallskip

In the next section, we prove of Theorem~1. Theorem~2 can be established analogously.

\section{Equal area and bounded perimeter\\---Proof of Theorem 1}
An unbounded convex planar region $C$ bounded by a segment and two {\em nonparallel} half-lines is said to be a {\em cone}. The segment is called the {\em base} of $C$, and the half-lines are called the {\em sides}. The \emph{angles} of $C$ are the angles between the base and the sides. The width of $C$ is defined as the minimum distance between two points on different sides. The width of $C$ is equal to the length of the base if and only if both angles of $C$ are at least $\pi/2$. Throughout this paper we consider closed regions, so we consider the triangles and the cones to be closed. In \emph{tiling} of the plane is the covering of the plane with pairwise internally disjoint closed regions.

We construct a tiling of the plane by unit area triangles, using the following infinite ``splitting procedure''. At each step of the procedure we obtain a tiling of the plane with finitely many triangles and cones.
\medskip

{\sc Splitting Procedure}

\smallskip

\noindent{\sc Step 0:} We split the plane by three half-lines starting from the origin $O$ into $3$ angular regions, the angle of each of which is $2\pi\over 3$. In each of these angular regions $A$, we choose a ``generic'' triangle $T=T(A)$ with unit area and all three angles at least $\pi/12$ such that one of its vertices is at the origin and the other two lie on opposite sides of $A$. The closure of $A\setminus T(A)$ is a cone. The {\em initial tiling} is the tiling of the plane consisting of the triangles and the cones constructed above: one triangle and one cone in each of the three angular regions.  (The precise definition of ``generic choice'' will be crucially important for the proof, and will be discussed later.)
\smallskip

\noindent{\sc Step $i+1$}, for $i\ge 0$: Let $C$ be one of the finitely many convex cones that we have after {\sc Step $i$}, which is closest to $O$.

{\bf If} {\em the width of $C$ is larger than $4$,} {\bf then} {\em split $C$ into two cones} with widths larger than $2$ by a half-line starting at a point $x$ of the base of $C$. It is clear that we can do this, further we have $2$ degrees of freedom to choose $x$ and the direction of the half-line emanating from $x$. We make a generic choice (see later).

{\bf If} {\em the width of $C$ is at most $4$,} {\bf then} there is a unit area triangle $T$ such that one of its sides coincides with the base of $C$ and its third vertex lies on a side of $C$ whose angle with the base is at most as large as the other angle of $C$. (We break ties arbitrarily.)  We {\em replace $C$ by the triangle $T$ and the cone $C'$ that is the closure of $C\setminus T$}.

\medskip

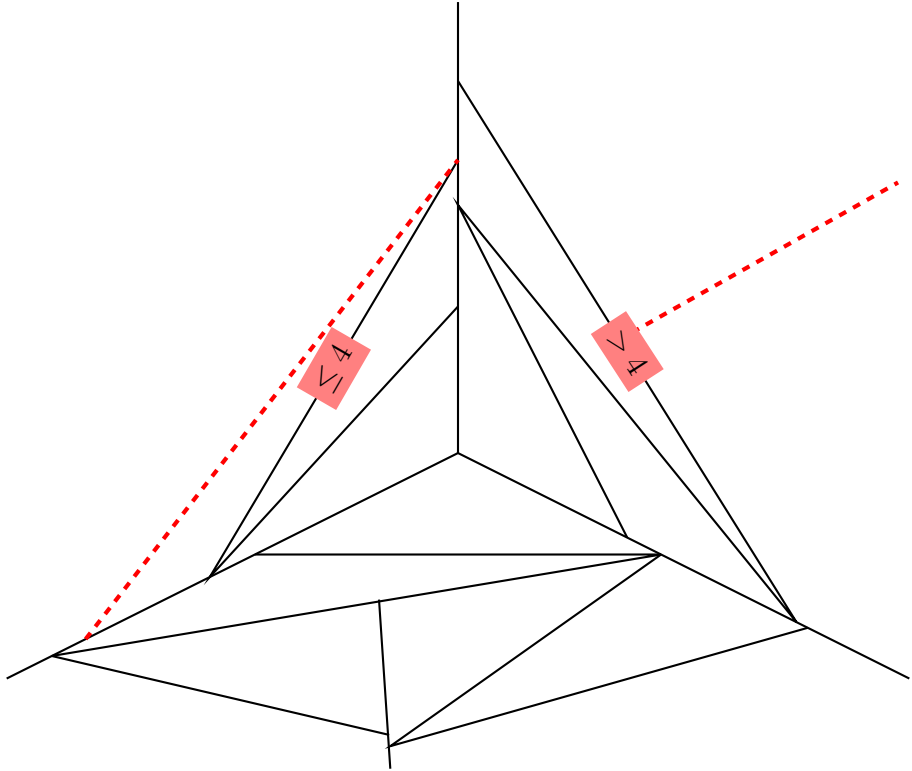
\begin{figure}\centering
\begin{tikzpicture}[scale=1.5]


\draw[thick] (-4,-2)-- (0,0) -- (0,4);
\draw[thick] (4,-2) -- (0,0);
\draw[thick] (0,2.6)--(-2.2,-1.1) -- (0,1.3);
\draw[thick] (1.5,-0.75) -- (0,2.2)-- (3,-1.5)--(0,3.3) ;
\draw[thick] (-1.8,-0.9) -- (1.8,-0.9)--(-3.6,-1.8);
\draw[ultra thick,dashed,red] (1.42,1) -- (3.9,2.4);
\draw[ultra thick,dashed,red] (-3.3,-1.65) -- (0,2.6);
\node[fill=Red!50, rotate=-57]  at (1.5, 0.9) {$>4$};
\node[fill=Red!50, rotate=60]  at (-1.1, 0.75) {$\le 4$};
\draw[thick] (-0.7,-1.3) -- (-0.6,-2.8);
\draw[thick] (1.8,-0.9)--  (-0.6,-2.6)--(3.1,-1.55);
\draw[thick] (-3.6,-1.8)--  (-0.61,-2.5);
\end{tikzpicture}
\caption{The {\sc Splitting Procedure}.}
\end{figure}
\medskip

\noindent{\bf Lemma 1.} {\em All cones $C$ created during the {\sc Splitting Procedure} satisfy
\smallskip

 {\bf Property P:} {The width of $C$  is larger than $2$, both of its angles are at most $11\pi/12$, and the sum of these two angles is at most $5\pi/3$.}}
\medskip

\noindent {\bf Proof.} Consider a cone $C$ of the initial tiling. It is the closure of $A\setminus T(A)$ for the angular region $A$. The angles of $C$ are $\pi$ minus an angle of $T(A)$, so they are at least $11\pi/12$. The sum of the angles of $C$ is $\pi$ plus the angle between the sides of $C$, and so it is exactly $5\pi/3$. Finally, the width of $C$ is the base of $C$, which is the longest edge of the unit area obtuse triangle $T(A)$, so it is longer than $2$. This establishes Property P for the cones in the initial tiling.

Our goal is to show that Property P is preserved throughout the whole {\sc Splitting Procedure}. This is obviously true for steps in which a cone $C$ of large width is split into two cones. Each of the resulting cones inherits one of its angles from $C$, while its other angle is strictly smaller than the other angle of $C$. The widths of the two new cones are required to be larger than $2$, so if the original cone had Property P, so do both of it parts.

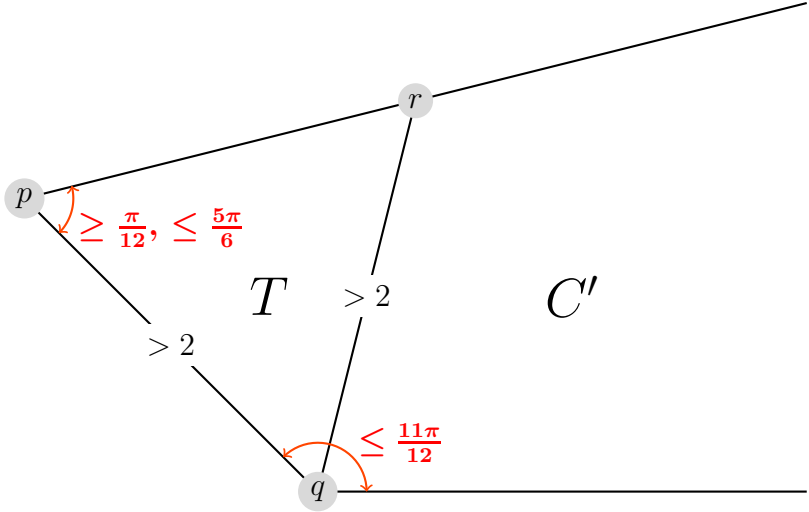
\begin{figure}\centering
\begin{tikzpicture}[scale=1.3]


\draw[thick] (5,5)--(-3,3)-- (0,0) -- (5,0);
\draw[thick] (1,4) -- (0,0);
\node[fill=white]  at (0.5, 2) {$>2$};
\node[fill=white]  at (-1.5, 1.5) {$>2$};
\node[circle,fill=Grey!30,inner sep = 2pt, minimum size=0pt]  at (0, 0) {$q$};
\node[circle,fill=Grey!30,inner sep = 2pt, minimum size=0pt]  at (-3, 3) {$p$};
\node[circle,fill=Grey!30,inner sep = 2pt, minimum size=0pt]  at (1, 4) {$r$};
\draw [<->, OrangeRed,thick,domain=-45:15] plot ({0.5*cos(\x)-3}, {0.5*sin(\x)+3});
\draw [<->, OrangeRed,thick,domain=0:135] plot ({0.5*cos(\x)}, {0.5*sin(\x)});
\node[circle,fill=white, inner sep = -22pt, minimum size=0pt]  at (0.85,0.5) {\textcolor{Red}{\bm{$\le\frac{11\pi}{12}$}}};
\node[circle,fill=white, inner sep = -22pt, minimum size=0pt]  at (-1.6,2.7) { \bm{\textcolor{Red}{$\ge\frac{\pi}{12},\, \le\frac{5\pi}{6}$}}};

\node[circle,fill=white, inner sep = -22pt, minimum size=0pt]  at (2.6,2) { \LARGE{$C'$}};

\node[circle,fill=white, inner sep = -22pt, minimum size=0pt]  at (-0.5,2) {\LARGE{$T$}};

\end{tikzpicture}
\caption{Illustration for the proof of Lemma 1. $\sphericalangle prq\ge\pi/12$.}
\end{figure}

It remains to verify that when a triangle $T$ is split off from a cone $C$ satisfying Property P, then the remaining cone $C'$ (the closure of $C\setminus T$) also satisfies Property P. Let us assume that $T=pqr$, where $pq$ is the base of $C$ and $r$ is on the side starting at $p$; see Fig. 2. The base of $C'$ is $qr$ and the angle of $C'$ at $q$ is smaller than the angle of $C$ at $q$, which was at most $11\pi/12$. The sum of the angles of $C'$ is the same as the sum of the angles of $C$, so at most $5\pi/3$. In view of the fact that the sides of $C'$ are contained in the sides of $C$, the width of $C'$ is at least as large as the width of $C$, which is larger than $2$. It remains to prove that the angle of $C'$ at $r$ is at most $11\pi/12$, or, equivalently, that $\sphericalangle prq$ is at least $\pi/12$.

Using the fact that the angles of $C$ add up to at most $5\pi/3$ and the angle of $C$ at $p$ is not larger than the other angle of $C$, we have $\sphericalangle qpr\le 5\pi/6$. As the other angle of $C$ is at most $11\pi/12$ and the sum is over $\pi$, we also have that $\sphericalangle qpr \ge \pi/12$. Among all unit area triangles $pqr$ with given base $pq$ and $\sphericalangle qpr\in [\pi/12, 5\pi/6]$,  $\sphericalangle prq$ attains the minimum for one of the extreme cases $\sphericalangle qpr=\pi/12$ or $\sphericalangle qpr=5\pi/6$.  If $pq$ had length $2$, both extremal configurations would correspond to unit area isosceles triangles with $\sphericalangle prq=\pi/12$. The length of $pq$ is at least as large as the width of $C$, which is larger than $2$. Therefore, $\sphericalangle prq\ge \pi/12$, as claimed.
$\Box$

\medskip

\noindent{\bf Lemma 2.} {\em All triangles created by the {\sc Splitting Procedure} have unit area and bounded perimeter. Altogether they form a tiling of the plane.}
\medskip

\noindent{\em Proof.} The unit area condition holds by the choice of these rectangles. For the bounded perimeter condition, we can safely ignore the three triangles in the initial tiling. So let $C$ be a cone with base $pq$, and let $T=pqr$ be the triangle cut off from $C$ by the {\sc Splitting Procedure} as in the end of the previous proof. The width of $C$ cannot exceed $4$, otherwise, $T$ would not have been cut off during the process. Notice that this width is at least the length of the base multiplied by $\sin(\sphericalangle qpr)\ge\sin(\pi/12)$, since $\sphericalangle qpr$ is the smaller angle of $C$ and it is still at least $\pi/12$. Therefore, the length of $pq$ is less than $16$, which together with $\sphericalangle prq\ge\pi/12$ (established in the previous proof) implies that the perimeter of $T$ is bounded.

The triangles created during the procedure are pairwise internally disjoint. We only have to verify that their union is the entire plane. Assume, for contradiction, that there is a point $x$ in the plane that belongs to a cone during the entire process. This means that every cone $C$ split during the process is at distance at most $xO$ from $O$. The base of $C$ is either a side of a triangle or is obtained by splitting the base of a previously created other cone. In either case, its length is bounded. Hence, there is a disk $D$ centered at $O$ which contains every triangle and the base of every cone created during the {\sc Splitting Procedure}. Take a larger disk $D'$ around $O$. As the angles of the cones are bounded away from $0$, there exists $\delta>0$ such that every cone covers a part of $D'$ whose area is at least $\delta$. Thus, at each stage of the procedure, the total number of triangles and cones is bounded from above by an absolute constant, which is a contradiction. $\Box$
\medskip

In the sequel, we spell out what we mean by ``generic choice'' during the {\sc Splitting Procedure}, and how this can be explored to avoid creating two congruent triangles.
\smallskip

There are precisely two ways to split off a triangle of unit area from a cone $C$ such that the remaining region is also a cone, depending on the side of $C$ which contains an edge of the triangle. Denote the resulting triangle and cone by $C^*_0$ and $C_0$, respectively, if the ``left'' side of $C$ contains an edge of the triangle, and by $C^*_1$ and $C_1$, respectively, if the ``right'' side. Extend this notation recursively to an arbitrary 0-1 sequence $s$ by letting
$$C_{s0}=(C_s)_0,\; \; \; C_{s1}=(C_s)_1,\; \; \;  C^*_{s0}=(C_s)^*_0,\; \; \; C^*_{s1}=(C_s)^*_1.$$
We call every triangle $C^*_s$ for some 0-1 sequence $s$ a {\em potential triangle} of the cone $C$. Note that if a cone $C$ appears in our process and its width is at most $4$, then in a future step one of the triangles $C^*_0$ or $C^*_1$ will be split off from $C$. If the width of the remaining cone still does not exceed $4$, we split off a further potential triangle of $C$, etc.

We need two simple facts, the first of which follows by trigonometric considerations, while the second is an elementary statement from real analysis.
\smallskip

\noindent{\bf Fact I.} {\em For any 0-1 sequence $s$, the length of each side of the triangle $C^*_s$ is an analytic function of the three parameters describing $C$: the length of its base and its two angles.}
\smallskip

\noindent{\bf Fact II.} {\em Let $\mathcal F$ be a countable collection of analytic functions $f$ defined on a connected open region $R$ of a Euclidean space. If none of these functions is identically zero, then every nonempty open subset of $R$ has a point where none of the functions vanishes.}

\smallskip

We would like to use our two degrees of freedom in choosing the splitting rays to guarantee that at every stage of the procedure all triangles that have already been created and all potential triangles of the existing cones are pairwise {\em noncongruent}. It will be sufficient to establish a slightly relaxed version of this property for the initial tiling: we allow potential triangles of the same cone to be congruent if they cannot possibly show up in the same tiling later, because they correspond to inconsistent future paths of splittings. To make this statement precise, we introduce the following definition.

\medskip
\noindent{\bf Definition.} A tiling of the plane into triangles and cones satisfies the {\em exclusion property} if there are no two congruent triangles among the triangles in the tiling and the potential triangles of the cones in the tiling, except possibly for pairs of potential triangles $C^*_s$ and $C^*_t$ obtained from the same cone $C$, and if two such triangles are congruent, then neither $s$ is an initial segment of $t$, nor vice versa.
\medskip

\noindent{\bf Lemma 3.} {\em For a suitable choice of the three triangles in {\sc Step 0}, the initial tiling satisfies the exclusion property.}
\medskip

\noindent{\em Proof.} Let $A$ be one of the angular regions in the initial tiling that will be split into a triangle $T$ and the cone $C$ that is the closure of $A\setminus T$. Let us write $A_s=C^*_s$ for a nonempty 0-1 string $s$, and let $A_s=T$ for the empty string $s$. We have one degree of freedom in choosing the length $l_A$ of the edge of $T$ on the left-hand side of $A$. The lengths of the other edges of $T$ and both angles of $C$ are analytic functions of $l_A$ and, by Fact~I, the same holds for the edge lengths of all potential triangles $A_s$. Hence, the lengths of all edges of these triangles in all three angular regions are analytic functions of the three free parameters $l_A$, $l_{A'}$, and $l_{A''}$ for the angular regions $A$, $A'$, and $A''$. Consider all pairwise differences between these functions. By Fact~II above, we can choose suitable values $l_A$, $l_{A'}$, and $l_{A''}$ within the required range in such a way that only those differences will vanish that happen to be identically zero for all possible choices of the parameters. We call such a choice of the parameters ``generic'', and argue that for a generic choice, the exclusion property holds. (Note that for a {\em generic} choice of the parameters, two potential triangles are congruent if and only if they are {\em identically} congruent, that is, they are congruent for {\em every} choice of the parameters.)

Observe that if $s$ is a proper initial segment of $t$, then the diameter of $A_s$ is strictly smaller than the diameter of $A_t$, so $A_s$ and $A_t$ cannot be congruent. Therefore, the exclusion property can only be violated by triangles belonging to distinct angular regions in the initial tiling. It is, therefore, sufficient to show that for any choice of distinct angular regions $A$ and $B$ in the initial tiling, and for every fixed pair of strings $s$ and $t$, the potential triangles $A_s$ and $B_t$ are not identically congruent.
For a fixed triangle of the initial tiling in $A$, the triangle $A_s$ is uniquely determined. Now fix the triangle of the initial tiling in $B$ in such a way that its diameter is larger than that of $A_s$. Then the diameter of $B_t$ larger than that of $A_s$, so they cannot be congruent. This shows that $A_s$ and $B_t$ are not identically congruent. $\Box$
\medskip

\noindent{\bf Lemma 4.} {\em If the exclusion property is satisfied for a finite tiling, then we can perform the next step in the {\sc Splitting Procedure} so that it is also satisfied for the resulting tiling.}
\medskip

\noindent{\em Proof.} If in the given step we split off a triangle $T$ from a cone $C$ of the current tiling, leaving another cone $D=C\setminus T$, then $T=C^*_i$ for $i=0$ or $i=1$ and $D^*_s=C^*_{is}$ for all 0-1 sequences $s$. Therefore, the exclusion property of the new tiling is inherited from the exclusion property of the old one.

If in the given step we split the cone $C$ into two cones, $D$ and $E$, then the resulting tiling does not necessarily have the exclusion property, because the potential triangles in $D$ and $E$ were not potential triangles in the old tiling. There are three types of possible conflicts, plus two more symmetric ones with the roles of $D$ and $E$ reversed (type (iii) is self-symmetric).

\begin{description}
\item[(i)] $D^*_s$ is congruent to a triangle $T$ that either appears in the current tiling or is a potential triangle of a cone other than $C$ in the current tiling.
\item[(ii)] $D^*_s$ and $D^*_t$ are congruent for some distinct 0-1 sequences $s$ and $t$, where $s$ is an initial segment of $t$.
\item[(iii)] $D^*_s$ and $E^*_t$ are congruent for some 0-1 sequences $s$ and $t$.
\end{description}

The requirement that the widths of the cones $D$ and $E$ are larger than $2$, determine an open region in the range of the two parameters that describe the ray separating $D$ and $E$. In view of Facts~I and II, it is sufficient to verify that any fixed conflict can be avoided for {\em some} choice of the ray. This implies that there is a suitable choice of the splitting line for which all possible conflicts will be avoided. According to Fact~I, the length of every edge of every relevant triangle is an analytic function of two free parameters that determine the split. As in the proof of Lemma~3, Fact~II guarantees that for a proper ``generic'' choice of the parameters, all coincidences between these lengths can be avoided that do not hold for every choice of the parameters, even if those choices for which two lengths differ do not belong to the permissible range of parameters (for example, if they give rise to cones of width smaller than $2$). For such a generic choice of the parameters, only those pairs triangles will be congruent that are identically congruent. It is enough to show that none of the conflicts listed above are unavoidable.

We will show that every conflict can be avoided by considering an ``extreme'' choice of the splitting line, very close to a side of $C$. We make the length of the base of $E$ and the angle between its sides both tend to zero. Then the parameters of $D$ will tend to those of $C$, hence the side lengths of the triangles $D^*_s$ will tend to the side lengths of $C^*_s$.

Therefore, conflicts of types (i) and (ii) cannot occur in the limit, as we assumed the tiling before the current step had the exclusion property. This implies that all conflicts of types (i) and (ii) can be avoided.

We have to be a bit more careful with conflicts of type (iii). Notice that the analytic dependence of the side lengths of $E_t$ remains valid even in the extreme case where $E$ has two parallel sides (in which case it is not a proper cone, but the definition extends to such shapes, too). Making the base length of $E$ go to zero, the width of $E$ will also tend to zero. Therefore, $E^*_t$ will degenerate in the limit, while $D^*_s$ will tend to $C^*_s$, which is a fixed nondegenerate triangle. Thus, this type of conflict can also be avoided. $\Box$
\smallskip

\medskip

According to Lemma~2, the {\sc Splitting Procedure} always yields a tiling of the plane with unit area and bounded perimeter triangles. Lemmas~3 and 4 show that with appropriate (generic) choices made during the procedure, all intermediate finite tilings will have the exclusion property. In particular, the resulting infinite tiling has no two congruent triangles. This completes the proof of Theorem~1.  $\Box\Box$


\begin{thebibliography}{DSS}
\bibitem{BBS10} I. B\'ar\'any, P. Blagojevi\'c, and A. Sz\H ucs: Equipartitioning by a convex 3-fan, {\em Advances in Math.} {\bf 223} (2010), 579--593.
\bibitem{BlZ14}  P. Blagojevi\'c and G. Ziegler: Convex equipartitions via equivariant obstruction theory, {\em Israel J. Math.} {\bf 200} (2014), no. 1, 49--77.
\bibitem{Fr16} D. Frettl\"oh: Noncongruent equidissections of the plane, {\tt  arXiv:1603.09132}.
\bibitem{KaHA14} R. Karasev, A. Hubard, and B. Aronov: Convex equipartitions: the spicy chicken theorem, {\em Geom. Dedicata} {\bf 170} (2014), 263--279.
\bibitem{KPT17} A. Kupavskii, J. Pach, and G. Tardos: Tilings with noncongruent triangles, preprint, {\tt arXiv:1711.04504}.
\bibitem{Na06} R. Nandakumar, Fair partitions. Blog entry, {\tt http://nandacumar.blogspot.de}, September 28, 2006.
\bibitem{Na14} R. Nandakumar: Filling the plane with non-congruent pieces, blog entries, December 2014, January 2015, June 2016, {\tt http://nandacumar.blogspot.in}.
\bibitem{NaR12} R. Nandakumar and N. Ramana Rao: Fair partitions of polygons: An elementary introduction, {\em Proc. Indian Academy of Sciences (Mathematical Sciences)} {\bf 122} (2012), 459--467.
\bibitem{So17} P. Sober\'on: Gerrymandering, sandwiches, and topology, {\em Notices Amer. Math. Soc.} {\bf 64} (2017), no. 9, 1010--1013.
\bibitem{Zi15} G. M. Ziegler: Cannons at sparrows, {\em Eur. Math. Soc. Newsl.} No. {\bf 95} (2015), 25--31.
\end{thebibliography}
\end{document}